\documentclass[a4paper,11pt]{amsart}
\usepackage[pdfdisplaydoctitle=true,
            colorlinks=true,
            urlcolor=blue,
            citecolor=blue,
            linkcolor=blue,
            pdfstartview=FitH,
            pdfpagemode=None,
            bookmarksnumbered=true
            ]{hyperref}
\usepackage{amssymb}
\usepackage{color}

\newtheorem{thm}{Theorem}[section]
\newtheorem{cor}[thm]{Corollary}
\newtheorem{lem}[thm]{Lemma}
\newtheorem{prop}[thm]{Proposition}

\theoremstyle{definition}

\theoremstyle{remark}
\newtheorem{rem}[thm]{Remark}
\newtheorem*{exple}{Example}

\newcommand{\R}{\mathbb{R}}
\newcommand{\Z}{\mathbb{Z}}
\newcommand{\N}{\mathbb{N}}
\newcommand{\C}{\mathbb{C}}
\newcommand{\e}{\mathbf{e}}

\newcommand{\Circle}{\mathbb{S}}

\newcommand{\id}{\mathrm{id}}

\newcommand{\g}{\mathfrak{g}}
\newcommand{\gstar}{\mathfrak{g}^{*}}

\newcommand{\DiffS}{\mathrm{Diff}^{\infty}(\mathbb{S})}
\newcommand{\DiffSstar}{\mathrm{Diff}_{1}^{\infty}(\mathbb{S})}

\newcommand{\CS}{\mathrm{C}^{\infty}(\mathbb{S})}
\newcommand{\CSzero}{\mathrm{C}^{\infty}_{0}(\mathbb{S})}

\newcommand{\D}[1]{\mathcal{D}^{#1}(\mathbb{S})}
\newcommand{\Dstar}[1]{\mathcal{D}_{1}^{#1}(\mathbb{S})}

\newcommand{\HH}[1]{H^{#1}(\mathbb{S})}
\newcommand{\HHzero}[1]{H^{#1}_{0}(\mathbb{S})}
\newcommand{\HHhat}[1]{\hat{H}^{#1}_{0}(\mathbb{S})}

\newcommand{\norm}[1]{\left\Vert#1\right\Vert}
\newcommand{\abs}[1]{\left\vert#1\right\vert}
\newcommand{\set}[1]{\left\{#1\right\}}

\DeclareMathOperator{\im}{Im} %

\hypersetup{
    pdftitle={The geometry of a vorticity model equation}, %  Title
    pdfauthor={J. Escher and B. Kolev and M. Wunsch}, % Author
    pdfsubject={MSC 2000: 35Q53, 58D05}, % Subject
    pdfkeywords={CLM equation, Euler equation, diffeomorphisms group of the circle}, % Keywords
    pdflang=EN % Language
    }

\begin{document}

\title[The geometry of a vorticity model equation]{The geometry of a vorticity model equation}

\author{Joachim Escher}
\address{Institute for Applied Mathematics, University of Hanover, D-30167 Hanover, Germany}
\email{escher@ifam.uni-hannover.de}
\author{Boris Kolev}
\address{LATP, CNRS \& Universit\'{e} de Provence, 39 Rue F. Joliot-Curie, 13453 Marseille Cedex 13, France}
\email{kolev@cmi.univ-mrs.fr}

%    Information for third author
\author{Marcus Wunsch}
\address{RIMS, Kyoto University, Kyoto 606-8502 Sakyoku Kitashirakawa Oiwakecho, Japan}
\email{mwunsch@kurims.kyoto-u.ac.jp}

\subjclass[2000]{58D05, 58B25, 35Q35, 37K65}
\keywords{Constantin-Lax-Majda equation, Euler equation on diffeomorphisms group of the circle} %

\date{October 22, 2010}%
%\dedicatory{}%
%\commby{}%

% ----------------------------------------------------------------

\begin{abstract}
We provide rigorous evidence of the fact that the modified Constantin-Lax-Majda equation modeling vortex and quasi-geostrophic dynamics \cite{OSW08} describes the geodesic flow on the subgroup $\DiffSstar$ of orientation-preserving diffeomorphisms $\varphi \in \DiffS$ such that $\varphi(1) = 1$ with respect to right-invariant metric induced by the homogeneous Sobolev space $\dot H^{1/2}(\Circle)$ and show the local existence of the geodesics in the extended group of diffeomorphisms of Sobolev class $H^{k}$ with $k\ge 2$.
\end{abstract}

\maketitle

% ----------------------------------------------------------------

\section{Introduction}
\label{sec:introduction}

Despite their derivation by Leonhard Euler as early as in 1757 \cite{Eul57}, the eponymous equations governing inviscid fluid flow still pose highly challenging mathematical problems.
For example, our knowledge about the propagation of regularity of solutions to the three-dimensional Euler equations remains fragmentary.
One way of attacking this and related question is to resort to simpler, lower-dimensional, differential equations whose solutions share features with those observed for (ideal) fluids.
This approach was chosen by Constantin, Lax, and Majda \cite{CLM85} in their derivation of a model equation for the three-dimensional vorticity equation (the Euler equation for the vorticity $\omega = \mbox{curl}\;u$ of the velocity vector $u$):
\begin{equation}\label{eq:CLM}
    \partial_t\omega = H \omega \; \omega, \quad t>0,\;x \in \R;
\end{equation}
the original equation being
\begin{equation*}
    \partial_t \omega + (u\cdot \nabla )\; \omega = \mathfrak D(\omega)\cdot \omega,\quad t>0,\;x \in \R^3.
\end{equation*}
This reduction occurs if the convective derivative $\partial_t + u\cdot \nabla$ is replaced by the temporal derivative, and it can be justified by the identical properties, in one or three space dimensions, respectively, of the singular integral operators $\mathfrak D$, given by the Biot-Savart law, and $H$, the Hilbert transform \cite{Pan96}.

A tremendous drawback of the vorticity model \eqref{eq:CLM}, however -- despite its having solutions exhibiting finite-time blow-up -- is the paradox that its viscous extension, first explored by Schochet \cite{Sch86}, has solutions which can break down earlier than in the inviscid regime.\footnote{Other viscous extensions of \eqref{eq:CLM} involving fractional Laplacians were proposed in \cite{WVM99} and in \cite{Sak03,Sak03a}. While these models have some physically reasonable properties, they fail to capture fundamental features of the corresponding 3D Navier-Stokes equations (cf. \cite{Sak03}).} In an attempt to circumvent this anomaly, De Gregorio \cite{DeG90} proposed another model for the vorticity equation:
\begin{equation}
    \partial_t\omega + u \omega_x = H\omega\; \omega,\quad u_x = H\omega,\quad t>0,\;x \in \R.
\end{equation}
In this convective perturbation of the CLM equation, he chose the ''velocity'' $u$ to be the antiderivative of the Hilbert transform of the ``vorticity'' $\omega$. Numerical studies \cite{OSW08} lead to the conjecture that De Gregorio's model equation has global solutions -- a fact that has yet to be verified mathematically. Also in \cite{OSW08}, a more general model equation of hydrodynamic type was presented:
\begin{equation}\label{eq:gCLM}
    \partial_t \omega + \alpha\; u \omega_x = H\omega\; \omega,\quad u_x = H\omega;
\end{equation}
$\alpha$ being an arbitrary real parameter. Three cases of which have been studied before:
\begin{itemize}
  \item $\alpha=-1$ corresponds to the model for the quasi-geostrophic equations of \cite{CCF05,CCF06};
  \item $\alpha=0$ reduces the model equation to the CLM equation;
  \item $\alpha = 1$ becomes the equation proposed by De Gregorio.
\end{itemize}

Recently, it has been shown that for any $\alpha$ lying in the negative half-line, there are solutions which blow up in finite time \cite{CC10}. However, the description of the asymptotic behavior in the case $\alpha > 0$ remains open to further scrutiny.

Note that the time scaling $t\mapsto t/\alpha$ transforms (\ref{eq:gCLM}) into the following equation
\begin{equation}\label{eq:mCLM}
    \partial_t \omega + u \omega_x + a\; u_x \omega = 0,\quad \omega = H u_x,
\end{equation}
with $a:=-1/\alpha$. This equation is known as the \textit{modified} CLM equation, cf. \cite{OSW08} for a recent study of (\ref{eq:mCLM}).  It was observed by Wunsch \cite{Wun10} that
in the case $a=2$ equation
(\ref{eq:mCLM}) admits a geometric interpretation as the description of the geodesic flow of the homogeneous $\dot H^{1/2}(\Circle)$ right-invariant metric on the homogeneous space $\DiffS/Rot(\Circle)$ of orientation-preserving diffeomorphisms of the circle modulo the subgroup of rotations $Rot(\Circle)$, putting it ``midway'' between the Burgers equation (giving the geodesic flow of the $L^2(\Circle)$ metric) and the Hunter-Saxton equation (describing the geodesic flow of the homogeneous $\dot H^1$ metric).

In this work, we will show that the periodic modified CLM can be realized as the geodesic flow of a symmetric linear connection
on the subgroup $\DiffSstar$ of orientation-preserving diffeomorphisms $\varphi \in \DiffS$ such that $\varphi(1) = 1$, which is canonically diffeomorphic to the coset manifold $\DiffS/Rot(\Circle)$.
For the value $a=2$ this connection is compatible with a Riemannian metric. In fact in this case, the metric is induced by the inertia operator $\Lambda:=HD$ (with respect to the $L^{2}$ inner product on the tangent bundle).

Moreover, we present a thorough study of the regularity properties of the geodesic flow on a suitable Banach approximation of the Fr\'{e}chet manifold $\DiffSstar$. Introducing furthermore Lagrangian coordinates, these regularity results allow us to prove the well-posedness of the modified CLM equation on large phase spaces.
An immediate application of our main results Theorem~\ref{theo:main} and Corollary~\ref{cor:mainthm}, in combination with Lemma~\ref{lem:H1_continuity} gives the following conclusion:

\begin{thm}\label{thm:mainCP}
Let $a\in\R$ and $k\ge 2$ be given. Then there exist $\delta_k>0$ and $T_k>0$ such for each $\omega_0\in \HH{k-1}$ with spatial mean zero and $\Vert \omega_0\Vert_{H^{k-1}}<\delta_k$ there exists a unique solution
\begin{equation*}
    \omega\in C( (-T_k,T_k), \HH{k-1} )\cap C^1((-T_k,T_k), \HH{k-2} )
\end{equation*}
to the modified CLM equation \eqref{eq:mCLM} with initial condition $\omega(0)=\omega_0$.
\end{thm}

Let us briefly outline the plan of our paper. In Section~\ref{sec:settings} we recall the
construction of Euler-Poincar\'{e} equations on Lie groups and provide basic facts on
Fourier multipliers on $\Circle$, which is the class of inertia operators
we are interested in. In Section~\ref{sec:euler} we realize the modified CLM equation as
an Euler equation on the Fr\'{e}chet Lie group $\DiffSstar$. Introducing Lagrangian
coordinates it is possible to re-formulate the geodesic flow on $\DiffSstar$ as a system
a first order ordinary differential equations on the tangent bundle of the Banach manifold
$\Dstar{k}$, consisting in all orientation-preserving diffeomorphisms $\varphi$ of Sobolev class
$H^k$ with $k\ge 2$ such that $\varphi(1)=1$. In Section~\ref{sec:regularity} we study the regularity of the vector field induced by the
above mentioned dynamical system. This section contains also our main results.
The proof of several continuity
properties of the composition mapping and some technical estimates for operators on Sobolev spaces are postponed to the Appendix.

% ----------------------------------------------------------------

\section{Settings}
\label{sec:settings}

A \emph{right-invariant} Riemannian metric on a Lie group $G$ is defined by its value at the unit element, that is by an inner product on the Lie algebra $\g$ of the group. If this inner product is represented by an invertible operator $A : \g \to \gstar$, for historical reasons, going back to the work of Euler on the motion of the rigid body, this inner product is called the \emph{inertia operator}.
The \emph{Levi-Civita connection} of such a Riemannian metric is itself \emph{right-invariant} and given by
\begin{equation}\label{eq:connection}
    \nabla_{\xi_{u}} \, \xi_{v} = \frac{1}{2} [\xi_{u},\xi_{v}] + B(\xi_{u},\xi_{v}),
\end{equation}
where $\xi_{u}$ is the right-invariant vector field on $G$ generated by $u\in \g$ and $B$ is the right-invariant tensor field on $G$, generated by the bilinear operator
\begin{equation}\label{eq:Christofel}
    B(u,v) = \frac{1}{2}\Big[ (\mathrm{ad}_{u})^{*}(v) + (\mathrm{ad}_{v})^{*}(u)\Big]
\end{equation}
where $u,v \in \g$ and $(\mathrm{ad}_{u})^{*}$ is the \emph{adjoint} (relatively to the inertia operator $A$) of the natural action of the Lie algebra on itself given by
\begin{equation*}
    \mathrm{ad}_{u} : v \mapsto [u, v].
\end{equation*}

Given a smooth path $g(t)$ in $G$, we define its \emph{Eulerian velocity}, which lies in the Lie algebra $\g$, by
\begin{equation*}
    u(t) = R_{g^{-1}(t)}\dot{g}(t)
\end{equation*}
where $R_{g}$ stands for the right translation in $G$. It can then be shown (see \cite{EK09} for instance) that $g(t)$ is a \emph{geodesic} if and only if its Eulerian velocity $u$ satisfies the first order equation
\begin{equation}\label{eq:Euler0}
    u_{t} = - B(u,u).
\end{equation}
known as the \emph{Euler equation} induced by $A$.

It was noticed in \cite{EK09} that the concept of Euler equation does not necessarily require the linear connection $\nabla$ to be Riemannian --- there may not exist a Riemannian metric which is \emph{preserved} by this connection. We have therefore called such an equation a \emph{non-metric Euler equation}. With this extended definition, every quadratic evolution equation on $\g$
\begin{equation*}
    u_{t} = Q(u),
\end{equation*}
corresponds to the reduced geodesic equation (Euler equation) of a \emph{right-invariant symmetric linear connection} on $G$.

The theory of Euler equations on a homogeneous space $G/K$ has been developed in \cite{KM03} in the metric case and more generally for Hamiltonian systems. It corresponds to a special case of the Hamiltonian reduction with respect to the subgroup $K$ action. Let $A : \g \to \gstar$, be the inertia operator of a \emph{degenerate} symmetric bilinear form $<\cdot , \cdot>$ on $\g$, such that $\ker A = \mathfrak{k}$, the Lie algebra of $K$. If moreover, the inner product $<\cdot , \cdot>$ is $\mathrm{Ad}_{K}$-invariant, that is
\begin{equation*}
    <\mathrm{Ad}_{k}u, \mathrm{Ad}_{k}v> = <u,v>,
\end{equation*}
for all $k \in K$ and all $u,v \in \g$, then $A$ induces a right $G$-invariant pseudo-Riemannian metric on the space $G/K$ of right cosets ($Kg$, $g\in G$). In that case, the Euler equation, which corresponds to the inertia operator $A$ and describes the geodesic flow on the homogeneous space $G/K$, has the following Hamiltonian form (see \cite{KM03}): it is the quotient with respect to the $K$-action of the restriction to $L  = \im A \subset \gstar$ of the following Hamiltonian equation on $\gstar$
\begin{equation*}
    m_{t} = - \mathrm{ad}_{A^{-1}m}^{*}m
\end{equation*}
for $m \in L$.

However, it appears difficult to work easily with a \emph{contravariant} formulation of this equation similar to equation~\eqref{eq:Euler0} in this more general situation. Indeed, in that case, the Eulerian velocity is only defined \emph{up to a path} in $K$ (see \cite{TV10} for a recent survey on the subject). Moreover, it is not clear how this formalism can be generalized to \emph{non-metric} Euler equations.

Fortunately, in the case we consider in this paper, these difficulties can be avoided because of a prolific structure. More precisely, in the situation we consider, there exists a closed subgroup $H$ of $G$, such that the restriction to $H$ of the canonical right action of $G$ on $G/K$ is transitive and without fixed points. In that case,
\begin{equation*}
    \g = \mathfrak{k} \oplus \mathfrak{h},
\end{equation*}
where $\mathfrak{k}$ is the Lie algebra of $K$ and $\mathfrak{h}$ is the Lie algebra of $H$ and the study of a degenerate, $\mathrm{Ad}_{K}$-invariant inner product on $\g$ with kernel $\mathfrak{k}$ can be reduced to an Euler equation on the Lie group $H$, where $\mathfrak{h}^{*}$ has be identified with
\begin{equation*}
    \mathfrak{k}^{0} = \set{m \in \gstar ; \; m(u) = 0, \, \forall u \in \mathfrak{k}}.
\end{equation*}

\begin{exple}
Let $E(3)$ be the Lie group of direct euclidean motions in 3-space. The homogeneous space $E(3)/SO(3)\simeq \R^{3}$ satisfies the hypothesis of our framework: the subgroup of translations $T(3)\simeq \R^{3}$ acts transitively and without fixed points on the quotient space. Notice however, that in this particular example, the subgroup $T(3)$ is a \emph{normal} subgroup of $E(3)$, something we do not assume explicitly in our more general framework.
\end{exple}

\begin{rem}
We emphasize, that contrary to what one might expect at first glance, the system of free motions of a \emph{rod} (degenerate rigid body) does not enter into this framework. Indeed, the configuration space of a rigid rod can be realized as the homogenous space $SO(3)/SO(2)\simeq S^{2}$ which is not diffeomorphic to any Lie group, otherwise, its tangent bundle would be trivial, which is not the case. Nevertheless, the geometric framework  of $\DiffS/Rot(\Circle)$ suits perfectly well for the study of the Hunter-Saxton equation (see \cite{Len08} for instance) and other hydrodynamical models we shall consider in this article.
\end{rem}

Let $\DiffS$ be the \emph{Fr\'{e}chet Lie group} of smooth and orientation preserving diffeomorphisms of the unit circle $\Circle \simeq \R/\mathbb Z$ and $\DiffS / Rot(\Circle)$ be the homogeneous space of right cosets
\begin{equation*}
    [\varphi] : = Rot(\Circle)\,\varphi ,
\end{equation*}
where $Rot(\Circle)$ is the subgroup of Euclidean rotations $x\mapsto x+s \pmod \Z$, $s\in \R$. The canonical \emph{right} action of the group $\DiffS$ on itself commutes with the left action of $Rot(\Circle)$ on $\DiffS$ and induces a \emph{right} action of $\DiffS$ on the quotient space $\DiffS/Rot(\Circle)$.

The restriction of this action to the subgroup $\DiffSstar$ of diffeomorphisms $\varphi \in \DiffS$ such that $\varphi(1) = 1$ is \emph{transitive} and \emph{simple} (without fixed point). Therefore, the restriction of the projection map $\varphi \mapsto [\varphi]$ to $\DiffSstar$ defines a bijection between $\DiffSstar$ and $\DiffS / Rot(\Circle)$. The inverse map is given by
\begin{equation*}
    [\varphi] \mapsto \varphi \cdot \varphi (1)^{-1}.
\end{equation*}
Notice however that the restriction to $\DiffSstar$ of the projection map is {\em{not a}} \emph{group morphism}. Otherwise $Rot(\Circle)$ would be a normal subgroup of $\DiffS$, which is not the case: $\DiffS$ is \emph{simple}, it has no (non trivial) normal subgroup \cite{GR07}.

The \emph{Fr\'{e}chet manifold} structure on $\DiffSstar$ is obtained by the existence of the global chart
\begin{equation}\label{eq:smooth_chart}
    U := \set{id + u; \; u \in \CS;\; u_{x}> -1, u(0)=0},
\end{equation}
which is an open set in the closed hyperplane $id + \CSzero$, where $\CSzero$ is the closed linear subspace
\begin{equation*}
    \CSzero:= \{ u \in \CS \big| u(0)= 0 \}.
\end{equation*}

The \emph{Lie bracket} on the tangent space at the unit element, identified with $\CSzero$, is given by
\begin{equation*}
    [u,v] = u_{x}v - uv_{x}
\end{equation*}
and we have
\begin{equation*}
    \CS = \CSzero \oplus \R,
\end{equation*}
where $\R$ is the Lie algebra of $Rot(\Circle)$.

\begin{rem}
To summarize, there is a prolific structure in the special framework we consider in this paper and which simplifies our work. This structure is essentially due to the fact that there exists a smooth section (not a group morphism however)
\begin{equation*}
    \DiffS / Rot(\Circle) \to \DiffSstar
\end{equation*}
of the canonical projection map
\begin{equation*}
    \DiffS \to \DiffS / Rot(\Circle).
\end{equation*}
In particular the Fr\'{e}chet Lie group  $\DiffS$ is diffeomorphic to the product manifold $Rot(\Circle) \times \DiffSstar$. Notice that the particular choice of the fixed point in the definition of $\DiffSstar$ does not affect the general structure.
\end{rem}

Consider now a \emph{non-negative} bilinear form on $\CS$ which can be written as
\begin{equation*}
    \langle u,v\rangle = \int_{\Circle} Au\cdot v\, dx,
\end{equation*}
where $A:\CS \to \CS$ is a linear, continuous, $L^{2}$-symmetric operator.

\begin{prop}\label{prop:riemannian_structure}
Suppose that $A$ satisfies the following three conditions
\begin{enumerate}
  \item $\ker A = \R$,
  \item $\im A = \set{ m \in \CS ;\; \int_{\Circle} m(x) \, dx = 0 }$,
  \item $AR_{s} = R_{s}A$, for all rotations $R_{s}$.
\end{enumerate}
Then, $A$ induces a \emph{weak Riemannian metric} on the homogeneous space $\DiffS / Rot(\Circle)$, identified with $\DiffSstar$. The
operator
\begin{equation*}
    B(u,v) = \frac{1}{2}A^{-1}\Big[ 2 A(v)u_{x} + A(v)_{x}u + 2 A(u)v_{x} + A(u)_{x}v \Big]
\end{equation*}
is well-defined on $\CSzero$ and the associated symmetric, right-invariant, linear connection on $\DiffSstar$ is compatible with the metric. The corresponding Euler equation on $\CSzero$ is given by
\begin{equation}\label{eq:Euler_metric}
    u_{t} = -B(u,u) = -A^{-1}\Big[ 2A(u)u_{x} + A(u)_{x}u\Big].
\end{equation}
\end{prop}

\begin{rem}
Notice that condition (3) is equivalent to the property for the degenerate inner product on $\CS$ defined by $A$ to be $\mathrm{Ad}_{Rot(\Circle)}$-invariant.
\end{rem}

\begin{rem}
Observe that the topology induced by the pre-Hilbertian structure on each tangent space of the Fr\'{e}chet manifold $\DiffS / Rot(\Circle)$ is weaker than the usual Fr\'{e}chet topology. For this reason such a structure is called a \emph{weak Riemannian metric}. On a Fr\'{e}chet manifold, only covariant derivatives along curves are meaningful. As expounded in \cite{EK09}, the general expression of a right-invariant, covariant derivative of a vector field $\xi(t) = (\varphi(t), w(t)) \in \DiffSstar \times \CSzero$ along the curve $\varphi(t) \in \DiffSstar$ is given by
\begin{equation*}
    \frac{D \xi(t)}{D t} = \left( \varphi, w_{t} + \frac{1}{2}[u,w] + B(u,w)\right),
\end{equation*}
where $u = \varphi_t \circ \varphi^{-1}$ and $B$ is a symmetric bilinear operator on $\CSzero$. However, and contrary to the finite dimensional case, the existence of a symmetric, linear connection on a Fr\'{e}chet manifold, compatible with a weak Riemannian metric, that is
\begin{equation*}
    \frac{d}{dt} \langle \xi , \eta \rangle_{\varphi} = \langle \frac{D \xi}{D t} , \eta \rangle_{\varphi} +\langle \xi ,  \frac{D \eta}{D t}\rangle_{\varphi},
\end{equation*}
is far from being granted.
\end{rem}

\begin{proof}[Proof of proposition~\ref{prop:riemannian_structure}]
If conditions (1) and (3) of proposition~\ref{prop:riemannian_structure} are fulfilled, $A$ induces a pre-Hilbertian structure on each tangent space of the homogeneous space $\DiffS / Rot(\Circle)$, identified with $\DiffSstar$. This inner product is given by
\begin{equation}\label{eq:Riemannian_metric}
    \langle \eta, \xi \rangle_{\varphi} = \langle \eta\circ\varphi^{-1}, \xi\circ\varphi^{-1} \rangle_{e} = \int_\Circle \eta\cdot  A_{\varphi}\xi\cdot \varphi_{x} \,dx,
\end{equation}
where $\eta, \xi \in T_{\varphi} \DiffSstar$ and $A_{\varphi} = R_{\varphi} \circ A \circ R_{\varphi^{-1}} $. This family of pre-Hilbertian structures, indexed by $\varphi \in \DiffSstar$, is smooth because composition and inversion are smooth on the Fr\'{e}chet Lie group $\DiffSstar$.
This way we obtain a right-invariant, \emph{weak Riemannian} metric on $\DiffSstar$.

Formula \eqref{eq:Christofel} cannot be used directly to define a connection compatible with the metric because the adjoint operators $\mathrm{ad}_{u}^{t}$ (relatively to the pre-Hilbertian structure) are not well-defined. Indeed, given $u,v,w \in \CSzero$, we have
\begin{equation*}
    <v, \mathrm{ad}_{u} w> = \int_{\Circle} \Big( 2A(v)u_{x} + A(v)_{x}u \Big)w \, dx
\end{equation*}
so that $\mathrm{ad}_{u}^{t}$ is well-defined if and only if $2A(v)u_{x} + A(v)_{x}u$ belongs to $\im A$, the space of smooth functions of mean value zero. One can check that this not the case in general. However, the expression
\begin{equation*}
    \big( 2 A(v)u_{x} + A(v)_{x}u \big) + \big( 2 A(u)v_{x} + A(u)_{x}v \big)
\end{equation*}
has mean value zero, and belongs to $\im A$ (condition (2)), provided $A$ commutes with $D$. This is true by virtue of lemma~\ref{lem:FOp}, if $A$ commutes with all rotations (condition (3)). Therefore, one can define
\begin{equation*}
    B(u,v) = \frac{1}{2}A^{-1}\Big[ 2 A(v)u_{x} + A(v)_{x}u + 2 A(u)v_{x} + A(u)_{x}v \Big]
\end{equation*}
and check that the associated right-invariant, symmetric linear connection on $\DiffSstar$ is compatible with the metric.
\end{proof}

More generally, for each $a\in \R$, the equation
\begin{equation}\label{eq:Euler}
    u_{t} = -B(u,u) = -A^{-1}\Big[ aA(u)u_{x} + A(u)_{x}u\Big].
\end{equation}
is the (non-metric) Euler equation of a well-defined symmetric, right-invariant, linear connection on $\DiffSstar$.

The special case where $A$ is a \emph{differential} operator with constant coefficients has been extensively studied (see, e.g., \cite{CK02,CK03,EK09}).
In this paper, we will need to extend the theory when $A$ is a \emph{Fourier multiplier}.
For later reference, let us first give a useful characterization of Fourier multipliers. Here and in the following
we use the notation $\e_{n}(x) = \exp(2\pi i n x)$, $n\in\mathbb{Z}$, $x\in \Circle$.

\begin{lem}\label{lem:FOp}
Let $P$ a \emph{continuous} linear operator on the Fr\'{e}chet space $\CS$. Then the following three conditions are equivalent:
\begin{enumerate}
  \item $P$ commutes with all rotations $R_{s}$.
  \item $[P,D] = 0$, where $D = d/dx$.
  \item For each $n \in \N$, there is a $p(n)\in\mathbb{C}$ such that  $P\e_{n} = p(n)\e_{n}$.
\end{enumerate}
In that case, we say that $P$ is a \emph{Fourier multiplier}.
\end{lem}

Since every smooth function on the unit circle $\Circle$ can be represented by
its Fourier series, we get that
\begin{equation}\label{eq:Fourier_series}
    (Pu)(x) = \sum_{k\in\mathbb{Z}} p(k)\hat u(k) \e_{k}(x),
\end{equation}
for every Fourier multiplier $P$ and every $u\in \CS$, where
\begin{equation*}
    \hat u(k):= \int_{\Circle} u(x) \e_{-k}(x)\,dx,
\end{equation*}
stands for the $k$-th Fourier coefficients of $u$.
The sequence
$p : \Z\to\C$
is called the \emph{symbol} of $P$.

\begin{proof}
Given $s \in \R$ and $u \in \CS$, let $u_{s}(x) := u(x+s)$. If $P$ commutes with translations we have
\begin{equation*}
    (Pu)_{s}(x) = (Pu_{s})(x).
\end{equation*}
Taking the derivative of both
sides
of this equation with respect to $s$ at $0$ and using the continuity of $P$, we get $DPu = PDu$ which proves the implication $(1) \Rightarrow (2)$.

If $[P,D] = 0$, then both $P\e_{n}$ and $\e_{n}$ are solutions of the linear differential equation
\begin{equation*}
    u^{\prime} = (-2\pi i n) u
\end{equation*}
and
are therefore equal
up to a multiplicative constant $p(n)$. This proves that $(2) \Rightarrow (3)$.

If $P\e_{n} = p(n)\e_{n}$, for each $n \in \N$ and $P$ is continuous, then we have representation~\eqref{eq:Fourier_series}. Therefore
\begin{align*}
    (Pu)_{s}(x) & = \sum_{k\in\Z} p(k)\hat u(k) \e_{k} (x + s) \\
        & = \sum_{k\in\Z} p(k)\widehat{u_{s}}(k) \e_{k}(x) = (Pu_{s})(x),
\end{align*}
which proves that $(3) \Rightarrow (1)$.
\end{proof}

\begin{rem}
Notice that the space of Fourier multipliers is a \emph{commutative subalgebra} of the algebra of linear operators on $\CS$ which contains all linear differential operators with constant coefficients.
\end{rem}

A \emph{Fourier multiplier} $P$ with symbol $p$ is said to be of order $s\in\N$ if there exists a constant $C >0$ such that
\begin{equation*}
   \abs{p(m)} \le \,C \abs{m}^{s},
\end{equation*}
for every $m \ne 0$. In that case, for each $k \ge s$, the operator $P$ extends to a bounded linear operator from $\HH{k}$
into $H^{k-s}(\Circle)$. We express this fact by  the notation $P\in\mathcal{L}(\HH{k},H^{k-s}(\Circle))$.

% ----------------------------------------------------------------

\section{The modified CLM equation as an Euler equation}
\label{sec:euler}

The homogeneous $\dot H^{1/2}$ norm defined on $\CSzero$ is introduced by means of Fourier series. We let
\begin{equation*}
	\norm{u}_{1/2}^{2} : = \sum_{k\in\Z} \vert k\vert\abs{\hat{u}(k)}^{2},
\end{equation*}
for $u\in \CS$. The corresponding inner product on $\CSzero$ can be written as
\begin{equation*}
    \langle u,v\rangle_e = \int_\Circle u\Lambda v \;dx,
\end{equation*}
where $u,v \in \CSzero$ and
\begin{equation*}
	\Lambda := H\circ D : \CSzero \to \CSzero^{*}
\end{equation*}
is  the inertia operator. In this formula, $D = d/dx$ and $H$ is the \emph{Hilbert transform}, defined either as a Cauchy principal value, cf. \cite{Pan96}
\begin{equation*}
(Hu)(x) = (p.v.) \int_{\Circle} u(x - y) \cot(\frac{y}{2})\; dy,
\end{equation*}
or, equivalently, as the \emph{Fourier multiplier} with symbol $h(k) = - i \, \mathrm{sgn}(k)$,
\begin{equation}\label{eq:Hilbert_transform_Fourier}
	(H u) (x) := - i \sum_{k=-\infty}^{+\infty}\mathrm{sgn}(k)\hat{u}(k)\e_{k}(x).
\end{equation}
The convention $\mathrm{sgn}(0) = 0$ permits to extend $H$ on $\CS$. Notice that
\begin{equation*}
    H^{2} = -Id
\end{equation*}
on $\CSzero^\ast$ and that $H$ defines a complex structure on this space. Moreover, $H$ is an isometry for the $L^{2}$ inner product of function equivalence classes having zero mean value.

Since $[\Lambda, D] = 0$ and the inertia operator $\Lambda$ is an isomorphism from $\CSzero$ onto
\begin{equation*}
    \CSzero^{*} := \set{ u \in \CS ;\; \hat{u}(0) = 0 },
\end{equation*}
the space of smooth functions of mean value zero, the existence of a linear connection compatible with the metric is granted and the Euler equation is defined.

\begin{thm}\label{thm:omega}
Given $a\in\mathbb{R}$,
the modified CLM equation
\begin{equation}\label{eq:genEuler}
    \omega_t + u \; \omega_x + a u_x\; \omega = 0, \quad \omega = \Lambda u,
\end{equation}
describes the geodesic flow of a right-invariant symmetric linear connection on the Fr\'{e}chet Lie group $\DiffSstar$. If $a=2$, the geodesic flow is metric and corresponds to the right-invariant homogeneous $\dot H^{1/2}(\Circle^1)$ metric.
\end{thm}

\begin{proof}
It suffices to replace $A$ by the expression $\Lambda = H \circ D$ in formula \eqref{eq:Euler} and to
use the definition $\omega = \Lambda u$, to get the general assertion.
\end{proof}

\begin{rem}

(a) Note that (\ref{eq:genEuler}) is equivalent to the Euler equation
\begin{equation}\label{eq:CLM_Euler}
    u_t = -\Lambda^{-1} \left[ u \; (\Lambda u)_x + a \Lambda (u) \; u_x \right].
\end{equation}
To the best of our knowledge, theorem~\ref{thm:omega} is the first time the
model equation for the 2D quasi-geostrophic and the Birkhoff-Rott equations studied in
\cite{CCF05,CCF06} has been identified as a non-metric Euler equation on $\DiffSstar$.

(b) We recall that the connection is Riemannian if $a=2$. Moreover, it follows from \cite{EW10} that for
\begin{equation*}
    a \not\in \set{-\frac{5}{3},\,-\frac{5}{4}, -\frac{5}{7}, \frac{1}{2} }
\end{equation*}
there is no inertia operator of Fourier multiplier type such that
(\ref{eq:genEuler}) can be realized as the geodesic flow with respect to the corresponding metric. These results extend similar statements for the $b$-equation \cite{Kol09,ES10}.

(c) It is also possible to consider evolution equations on $\DiffSstar$, related to the inertia operator $\Lambda^{2} = -D^2$. Given $a\in\R$, one may study the family
\begin{equation*}
    m_t + um - a\,u_x m = 0.
\end{equation*}
The most prominent equations in this family are the Hunter-Saxton equation
\begin{equation}\label{eq:HS}
    m_t+um +2u_x m=0,
\end{equation}
and the Proudman-Johnson equation
\begin{equation}\label{eq:PJ}
    m_t+um-u_x m=0,
\end{equation}
respectively.
The Hunter-Saxton equation is closely related to the Camassa-Holm Equation and the geometric picture of (\ref{eq:HS}) is in fact fairly good understood, cf. \cite{Len07}. In particular, (\ref{eq:HS}) is the geodesic flow on $\DiffSstar$ with respect to the inertia operator $-D^2$, i.e. with respect to the homogeneous $\dot H^1$-metric on $\CSzero$. In contrast, (\ref{eq:PJ}) is a non-metric Euler equation, cf. \cite{EW10}.
\end{rem}

% ----------------------------------------------------------------

\section{The geodesic flow on $\Dstar{k}$}
\label{sec:regularity}

In this section, we will study the regularity of the geodesic flow on suitable \emph{Banach approximations} of the Fr\'{e}chet manifold $\DiffSstar$. More precisely, let $\D{k}$ be the Banach manifold of orientation-preserving diffeomorphisms $\varphi$ of Sobolev class $H^k$ (defined for some integer $k \ge 2$). This Banach manifold is a \emph{topological group} with respect to composition of diffeomorphisms but it is not a \emph{Lie group}. Indeed, on $\D{k}$, right translation $R_{\varphi} : \psi \mapsto \psi\circ\varphi$
is linear, hence
smooth; whereas left translation $L_{\varphi} : \psi \mapsto \varphi\circ\psi$ is only continuous but not  differentiable
(see \cite{EM70,Ham82}).

\begin{rem}
The space of homeomorphisms of the circle of Sobolev class $H^{1}$ (as well as their inverse) is not a group. Indeed, let $1/2 < \alpha < 1/\sqrt{2}$. Then $F: x \mapsto x^{\alpha}$ is an increasing homeomorphism of $[0,1]$ which induces an homeomorphism of the circle. One can check that $F$ as well as $F^{-1}$ are of class $H^{1}$ but that $F \circ F$ is not, since
$\partial(F\circ F)\not\in L^{2}$.
\end{rem}

Analogous to $\DiffSstar$, the codimension one Banach submanifold $\Dstar{k}$ of diffeomorphisms $\varphi \in \D{k}$ such that $\varphi(1) = 1$ is covered by the global chart
\begin{equation}\label{eq:Hk_chart}
    U_{k} = \set{id + u; \; u \in \HH{k}, u_{x}> -1, u(0)=0},
\end{equation}
which is an open set in the closed hyperplane $id + \HHzero{k}$, where $\HHzero{k}$ is the closed linear subspace
\begin{equation*}
   \HHzero{k} := \{ u \in \HH{k}  ;\;  u(0)= 0 \}.
\end{equation*}
Since the manifold $\Dstar{k}$ is described by a global chart which is an open set of the vector space $\HHzero{k}$, its fiber bundle is trivial and all tangent spaces $T_{\varphi}\Dstar{k}$ can be identified canonically with $\HHzero{k}$, using this chart.

As a Fourier multiplier of order 1, the inertia operator $\Lambda = H \circ D$ extends to a bounded linear isomorphism
\begin{equation*}
   \Lambda : \HH{k} \to \HH{k-1}.
\end{equation*}

Notice that the restriction of $D$ to $\HHzero{k}$ is a bounded isomorphism onto
\begin{equation*}
   \HHhat{k} := \set{ m \in \HH{k} ;\;  \int_{\Circle} m \, dx=0}.
\end{equation*}
Its inverse is given by
\begin{equation*}
   D^{-1} : m \mapsto u ;\; u(x) = \int_{0}^{x} m(t) \, dt ,\quad x\in \Circle.
\end{equation*}
Since the Hilbert transform $H$, which is an isometry for the $L^{2}$ product for zero-mean periodic functions, commutes with $D$, its restrictions to $\HHhat{k}$ is an isometry of $\HHhat{k}$ (for the $H^{k}$ inner product).

The Euler equation \eqref{eq:CLM_Euler} is not an ODE on $\HHzero{k}$ because the second-order term $\Lambda u_x=H u_{xx}$ is not regularized by the inverse of the first order Fourier multiplier $\Lambda$. By introducing Lagrangian coordinates, however, one can get around this impediment and it is possible to re-formulate \eqref{eq:CLM_Euler} as a well defined vector field on the Banach manifold $\Dstar{k}\times \HHzero{k}$.

\begin{thm}\label{theo:main}
Let $a\in\R$ and $k\in\mathbb{N}$ with $k\ge 2$ be given.
The time-dependent vector field $u\in \HHzero{k}$ is a solution to the modified CLM equation if and only if $(\varphi,v)$ is a solution to
\begin{equation}\label{eq:CLM_Cauchy}
    \begin{cases}
        \ \partial_t \; \varphi = v \\
        \ \partial_t \; v = S_\varphi(v),
    \end{cases}
\end{equation}
where $S_\varphi(v) := (R_\varphi \circ S \circ R_{\varphi^{-1}})(v)$,
and
\begin{equation*}
    S(u) = \Lambda^{-1} \left\{ [\Lambda,u] u_x - a (\Lambda u)u_x \right\}.
\end{equation*}
Moreover, the second order vector field
\begin{equation*}
    \Phi: \Dstar{k} \times \HHzero{k} \rightarrow \HHzero{k} \times \HHzero{k}
\end{equation*}
given by $\Phi(\varphi,w) = \left( w, S_\varphi(w)\right)$ is of class $C^{\infty}$.
\end{thm}

\begin{rem}
The second order vector field $\Phi$ is called the \emph{spray} (of the metric or the linear connection) in the literature.
\end{rem}

The first part of the theorem results from the following observation. Let $u$ be a time-dependent vector field on $J\times \Circle$, where $J$ is an open interval in $\R$, and let $\varphi$ be its flow, i.e. $\varphi_t = u \circ \varphi$. Setting $v = u \circ \varphi$, we get $v_t = (u_t + u u_x)\circ \varphi$ by the chain rule. Hence $u$ is a solution to \eqref{eq:CLM_Euler} if and only if
\begin{eqnarray*}
    u_t + u u_x &=& -\Lambda^{-1} \left[ u (\Lambda u)_x - \Lambda (u u_x) + a (\Lambda u) u_x \right]\\
    &=& \Lambda^{-1} \left\{ [\Lambda,u]u_x - a (\Lambda u)u_x\right\}.
\end{eqnarray*}

The proof of the second part of the theorem, i.e. the smoothness of the spray $\Phi$ consists of several reductions, some of them
being true for general Fourier multipliers. We outline these  reductions in the remainder of this section. Some technicalities will
be postponed to Appendix~\ref{sec:proof_prop}.

Before entering into the details of the proof, let us state that the above result allows us to apply the Picard-Lindel\"{o}f theorem, which immediately yields:

\begin{cor}\label{cor:mainthm}
Let $a\in\R$ and $k\ge 2$ be given. Then there exist $\delta_k>0$ and $T_k>0$ such for each $u_0\in \HHzero{k}$ with $\Vert u_0\Vert_{H^k}<\delta_k$ there exists a unique solution
\begin{equation*}
    (\varphi,v)\in C^{\infty}( (-T_k,T_k), \Dstar{k} \times \HHzero{k} )
\end{equation*}
to \eqref{eq:CLM_Cauchy} such
that
$\varphi(0) = \id$ and $v(0)=u_0$.
\end{cor}

Let us start with the first reduction. If we assume that the conjugation  $A_\varphi$ of the inertia operator $A$ is of class
$C^m$ then the spray $\Phi$ is of class $C^{m-1}$.

\begin{prop}\label{prop:redu1}
Let $m \ge 1$, $a \in \R$, $s \ge 1$ and $k\ge s +1$. Let $A$ be a Fourier multiplier of order $s$. Suppose that
\begin{equation*}
	(\varphi, v) \mapsto A_{\varphi} (v) = R_{\varphi} \circ A \circ R_{\varphi^{-1}} (v).
\end{equation*}
is of class $C^{m}$ from $\D{k} \times \HH{k}$ to $\HH{k-s}$ and that $A$ induces an isomorphism from $\HHzero{k}$ onto $\HHhat{k-s}$. Then
\begin{equation*}
	(\varphi, v) \mapsto S_{\varphi} (v) = R_{\varphi} \circ S \circ R_{\varphi^{-1}} (v)
\end{equation*}
where
\begin{equation*}
    S(u) = A^{-1} \left\{ [A,u] u_x - a (A u)u_x \right\},
\end{equation*}
is of class $C^{m-1}$ from $\Dstar{k} \times \HHzero{k}$ to $\HHzero{k}$.
\end{prop}
\begin{proof}
Let  $P(u) := (A u) u_x$ and $Q(u) := [A,u] u_x$.
We have
\begin{equation*}
    S_{\varphi}(v) = A_{\varphi}^{-1} \left\{ Q_{\varphi}(v) - aP_{\varphi}(v) \right\},
\end{equation*}
where the subscript $\varphi$ indicates the conjugacy by the right translation $R_{\varphi}$ in $\Dstar{k}$. Although $P$ and $Q$ are smooth operators, these results do not carry over when conjugated with translation in $\Dstar{k}$ since for $k\ge 2$ these sets only form topological groups: neither composition nor inversion are differentiable.

Given an operator $K$, we introduce the following notation
\begin{equation*}
	\tilde{K}(\varphi, v) := (\varphi, K_{\varphi}(v)),
\end{equation*}
where $K_{\varphi}(v) = R_{\varphi} \circ K \circ R_{\varphi^{-1}} (v)$.

1) We have $P_{\varphi}(v) = \big( A_{\varphi} (v) \big)\big(  D_{\varphi}(v) \big)$. But
\begin{equation*}
	(\varphi,v) \mapsto  D_{\varphi}(v)
\end{equation*}
is smooth since $D_{\varphi}(v) = v_{x}/\varphi_{x}$ and $\HH{k}$ is a Banach algebra for $k \ge 1$.
Also $H^{k-s}(\Circle)$ is a Banach algebra because $k-s \ge 1$. Hence the fact that $P_{\varphi}(v) \in \HH{k-s}$
and our assumption ensure that
\begin{equation*}
	\tilde{P} :  \D{k} \times \HH{k} \to \D{k} \times \HH{k-s},
\end{equation*}
is of class $C^{m}$.

2) Since
\begin{equation*}
    d_{(\varphi,v) }\tilde{A} (\delta\varphi, \delta v) = \left(
    \begin{array}{cc}
 			\id &  0\\
			 * & A_{\varphi}
	\end{array}
	\right)
\end{equation*}
is a bounded, linear, invertible operator from $\HHzero{k} \times \HHzero{k}$ to $\HHzero{k} \times \HHhat{k-s}$, we conclude, using the inverse mapping theorem on Banach spaces, that
\begin{equation*}
	\tilde{A} ^{-1}:  \Dstar{k} \times\HHhat{k-s} \to \Dstar{k} \times \HHzero{k}
\end{equation*}
is of class $C^{m}$.

3) Taking $P = A$ and $\delta\varphi_{1} = v = u \circ \varphi$ in Proposition~\ref{lem:nth_derivative} when $\varphi, v$ are smooth, we get
\begin{equation*}
	\partial_{\varphi}A_{\varphi} (\varphi,v,v) = \{ [u,A] \circ D \}_{\varphi}(u \circ \varphi) = - Q_{\varphi}(v).
\end{equation*}
Now since smooth maps are dense in Sobolev spaces, this relation is still valid for $\varphi \in \Dstar{k}$ and $v \in \HHzero{k}$ and therefore
\begin{equation*}
	\tilde{Q} :  \Dstar{k} \times \HHzero{k} \to \Dstar{k} \times \HHhat{k-s},
\end{equation*}
is of class $C^{m-1}$.
The assertion now follows from the chain rule.
\end{proof}

Next we show that the conjugation of an inertia operator of Fourier multiplier type is in fact smooth. In order to do so we first consider operators in the smooth category and extend them in s second step to Sobolev spaces.

Let $(\varphi,v) \mapsto P_{\varphi}(v)$ be a smooth map on the Fr\'{e}chet manifold $\DiffS \times \CS$, where $P$ is linear in $v$. The partial G\^{a}teaux derivative of $P$ in the first variable $\varphi$ and in the direction $\delta\varphi_{1}\in\CS$ is a smooth map which is linear both in $v$ and $\delta\varphi_{1}$ and that we will denote by
\begin{equation}\label{eq:Gateaux_derivative}
	\partial_{\varphi}P_{\varphi} (v,\delta\varphi_{1}).
\end{equation}
Therefore, the partial G\^{a}teaux derivative of $P$ in the variable $\varphi$ is a map of three independent variables : $\varphi$, $v$, $\delta\varphi_{1}$. The second partial derivative of $P$ is directions $\delta\varphi_{1}, \delta\varphi_{2} \in \CS$ is the partial G\^{a}teaux derivative of \eqref{eq:Gateaux_derivative} in the variable $\varphi$ and in the direction $\delta\varphi_{2}$. We will denoted it by
\begin{equation*}
	\partial^{2}_{\varphi}P_{\varphi} (v,\delta\varphi_{1},\delta\varphi_{2}).
\end{equation*}
It can be checked that this expression is symmetric in $\delta\varphi_{1},\delta\varphi_{2}$ (see \cite{Ham82}). Inductively, we define this way the n-th partial derivative of $P$ in directions $\delta\varphi_{1}, \dotsc , \delta\varphi_{n}$ and we write it as
\begin{equation*}
	\partial^{n}_{\varphi}P_{\varphi} (v,\delta\varphi_{1}, \dotsc , \delta\varphi_{n}).
\end{equation*}
The space of linear operators on a Fr\'{e}chet space is a locally convex topological vector space, but in general is
not a Fr\'{e}chet space (see \cite{Ham82}). For this reason, we will avoid taking limits and derivatives of linear operators. In the sequel, if such equalities appear for notational simplicity, it just means equality of operators.

\begin{prop}\label{lem:nth_derivative}
Let $P$ be a continuous, linear operator on $\CS$ and let
\begin{equation*}
	P_{\varphi} = R_{\varphi}PR_{\varphi}^{-1},
\end{equation*}
where $\varphi \in \DiffS$. Then, given $n\in\mathbb{N}$, we have
\begin{equation}\label{eq:nth_derivative}
    \partial^{n}_{\varphi}P_{\varphi} (v,\delta\varphi_{1}, \dotsc ,\delta\varphi_{n}) = R_{\varphi}P_{n}(u_{1}, \dotsc , u_{n})R_{\varphi}^{-1}(v),
\end{equation}
where $u_{i} = \delta\varphi_{i}\circ \varphi^{-1}$ and $P_{n}$ is the multilinear operator defined inductively by $P_{0} = P$ and
\begin{multline}\label{eq:recurrence_relation}
    P_{n+1}(u_{1}, \dotsc , u_{n+1}) = [u_{n+1}D,P_{n}(u_{1}, \dotsc , u_{n})] \\
    - \sum_{i=1}^{n}P_{n}(u_{1}, \dotsc ,u_{i,x}u_{n+1}, \dotsc , u_{n}).
\end{multline}
\end{prop}

\begin{rem}
For a Fourier multiplier, that is, if $[P,D] = 0$, we have
\begin{equation*}
    P_{1}(u_{1}) = [u_{1},P]D,
\end{equation*}
and
\begin{equation*}
    P_{2}(u_{1},u_{2}) = [u_{1},[u_{2},P]]D^{2} + [u_{1},P][u_{2},D]D + [u_{2},P][u_{1},D]D.
\end{equation*}
\end{rem}

\begin{proof}
Formula~\eqref{eq:nth_derivative} is trivially true for $n=0$. Now suppose it is true for some $n\in\N$, that is
\begin{equation*}
    \partial^{n}_{\varphi}P_{\varphi} (v,\delta\varphi_{1}, \dotsc ,\delta\varphi_{n}) = R_{\varphi}P_{n}(u_{1}, \dotsc , u_{n})R_{\varphi}^{-1}(v),
\end{equation*}
where $u_{i} = \delta\varphi_{i}\circ \varphi^{-1}$ for $1 \le i \le n$. Notice that, for fixed $\delta\varphi_{1}, \dotsc ,\delta\varphi_{n}$
\begin{equation*}
    P_{n}(u_{1}, \dotsc , u_{n}) = P_{n}(\delta\varphi_{1}\circ \varphi^{-1}, \dotsc , \delta\varphi_{n}\circ \varphi^{-1})
\end{equation*}
is a family of linear operator on $\CS$ indexed by $\varphi$ and which depend on $\varphi$ only through the $u_{i}$. Let $\varphi(s)$ be a smooth path in $\DiffS$ such that
\begin{equation*}
    \varphi(0) = \varphi, \qquad \partial_s \; \varphi(s) \big|_{s = 0} = \delta\varphi_{n+1}
\end{equation*}
and let $u_{n+1} = \delta\varphi_{n+1}\circ\varphi^{-1}$. We compute first
\begin{equation*}
	\dot{R}_{\varphi} := \partial_s \; R_{\varphi(s)} \big|_{s = 0} = R_{\varphi}u_{n+1}D,
\end{equation*}
so that
\begin{equation*}
	R_{\varphi}^{-1}\dot{R}_{\varphi} = u_{n+1}D,
\end{equation*}
and
\begin{equation*}
    \dot{u}_{i} := \partial_s \; \left(\delta\varphi_{i}\circ\varphi(s)^{-1}\right) \big|_{s = 0} = - u_{i,x}u_{n+1},
\end{equation*}
for $1\le i \le n$. We have then
\begin{equation*}
    \dot{P}_{n} := \partial_s \; P_{n}(u_{1}, \dotsc , u_{n}) \big|_{s = 0} = - \sum_{i=1}^{n}P_{n}(u_{1}, \dotsc ,u_{i,x}u_{n+1}, \dotsc , u_{n}).
\end{equation*}
Finally, we have (simplifying the notation $P_{n}$ for $P_{n}(u_{1}, \dotsc , u_{n})$)
\begin{align*}
    \partial_s \; R_{\varphi}P_{n}R_{\varphi}^{-1} \big|_{s = 0}
    & = \dot{R}_{\varphi}P_{n}R_{\varphi}^{-1} + R_{\varphi} \dot{P}_{n} R_{\varphi}^{-1} - R_{\varphi} P_{n} \left(R_{\varphi}^{-1} \dot{R}_{\varphi} R_{\varphi}^{-1}\right) \\
    & = R_{\varphi}\left( R_{\varphi}^{-1}\dot{R}_{\varphi}P_{n} - P_{n}R_{\varphi}^{-1}\dot{R}_{\varphi} \right)R_{\varphi}^{-1} + R_{\varphi} \dot{P}_{n} R_{\varphi}^{-1} \\
    & = R_{\varphi} \left( [u_{n+1}D , P_{n}] + \dot{P}_{n} \right)R_{\varphi}^{-1},
\end{align*}
which gives the recurrence relation~\eqref{eq:recurrence_relation}, since
\begin{equation*}
    \partial^{n+1}_{\varphi}P_{\varphi} (v,\delta\varphi_{1}, \dotsc ,\delta\varphi_{n+1}) = \partial_s \; \left(R_{\varphi}P_{n}(u_{1}, \dotsc , u_{n})R_{\varphi}^{-1}(v)\right) \big|_{s = 0},
\end{equation*}
the proof, the proof is complete.
\end{proof}

Proposition~\ref{lem:nth_derivative} is the core of the following result, which ensures smoothness of the inertia operator $\Lambda_\varphi(v)$ in both variables with respect to suitable Sobolev norms. To avoid too much technicalities here, we postpone its proof to Appendix B.

\begin{prop}{\em(Smoothness of the conjugate of the inertia operator)}\label{prop:smoothness_of_lambda}
Let $k\ge 2$ and $\Lambda = H \circ D$. Then
\begin{equation*}
	(\varphi, v) \mapsto \Lambda_{\varphi} (v) = R_{\varphi} \circ \Lambda \circ R_{\varphi^{-1}} (v).
\end{equation*}
is of class $C^{\infty}$ from $\D{k} \times \HH{k}$ to $\HH{k-1}$.
\end{prop}

In contrast to finite dimensional Riemannian geometry the topology of the fibre of the tangent bundle is fundamental importance
in the infinite dimensional case.
It is clear that in the smooth category the pre-Hilbertian structure defined by \eqref{eq:Riemannian_metric} will not induce the Fr\'{e}chet topology of the tangent space $\CSzero$. The very same is  true if we complete the tangent space with respect to a general Banach norm. Therefore we call the metric induced by \eqref{eq:Riemannian_metric} a {\em{weak Riemannian metric}}.

\begin{cor}{\em(Smoothness of the metric and the spray)}
For each $k\ge 2$, the right-invariant, weak Riemannian metric defined by formula~\eqref{eq:Riemannian_metric} on $\DiffSstar$ with $A=\Lambda$ extends to a \emph{smooth} weak Riemannian metric on the Banach manifold $\Dstar{k}$ with a smooth geodesic spray.
\end{cor}

\begin{rem}
To conclude this section, it could be worth to bring together the present work with the right-invariant metric defined by the inertia operator
\begin{equation*}
    A := HD(D^{2}-1)
\end{equation*}
defined on the diffeomorphism group of the circle which fixes the three points $-1,0,1$. This metric has been related with the Weil-Petersson metric on the universal Teichm\"{u}ller space $T(1)$ in \cite{TT06}. The corresponding geodesic flow has been extensively studied in \cite{GB09}. Recall first that $\D{s}$, the space of homeomorphisms of class $H^{s}$ as well as their inverse is a topological group only for $s>3/2$ and that $3/2$ is therefore a critical exponent. One of the main results in \cite{GB09} is that, the inertia operator $A$ defines on a suitable replacement for the ``$H^{3/2}$ diffeomorphism group", a right-invariant \emph{strong Riemannian structure} which is moreover \emph{complete} (geodesics are defined globally).

Our point of view in this paper is completely different in the sense that we work on a well defined topological group $\D{s}$ for $s >3/2$ equipped with a Banach manifold structure\footnote{We deliberately decided to restrict to $s\in \N$ for simplicity but choosing $s\in \R$ does not invalidate our results.}. The price to pay for this nice structure is the fact that the metric only defines a \emph{weak Riemannian structure}. Nevertheless, we have been able to show local existence of the geodesics, also in this context.
\end{rem}

% ----------------------------------------------------------------

\appendix

% ----------------------------------------------------------------

\section{Continuity lemmas}
\label{sec:continuity_lemmas}

In this section we provide some continuity properties of the composition mapping in  Sobolev spaces.
Given Fr\'{e}chet spaces $X$ and $Y$, let $\mathcal{L}(X,Y)$ denote the space of all continuous linear operators from $X$ into $Y$.

\begin{lem}\label{lem:sep_continuity}
Let $X$, $Y$ be Fr\'{e}chet spaces and let $G$ be a metric space. Given $F\,:\, G\times X\to Y$, assume that
\begin{equation}\label{eq:sep_cont1}
\begin{cases}
     F(g,\cdot)\in\mathcal{L}(X,Y)\quad\text{for all}\quad g\in G, \\
     F(\cdot,x)\in C(G,Y) \quad\text{for all}\quad x\in X.
\end{cases}
\end{equation}
Then $F\in C(G\times X,Y)$.
\end{lem}
\begin{proof}
Fix $(g_0,x_0)\in G\times X$ and pick a sequence $(g_n,x_n)$ in $G\times X$ such that $\lim_n(g_n,x_n)=(g_0,x_0)$. Let further $V$ denote a neighbourhood of $F(g_0,x_0)$ in $Y$. We set
\begin{equation*}
    B_n:= F(g_n,\cdot)\in\mathcal{L}(X,Y),\quad n\in\N.
\end{equation*}
Then, given $x\in X$, we have
\begin{equation*}
    \lim_n B_n(x)=\lim_n F(g_n,x)=F(g_0,x).
\end{equation*}
Hence $\{B_n(x)\,;\, n\in\N\}$ is bounded in $Y$. Invoking the \emph{uniform boundedness principle} in Fr\'{e}chet spaces (see \cite[Theorem II.11]{DS58}), we deduce that the family $\{B_n\,;\, n\in\N\}$ is equicontinuous. In particular there is a neighbourhood $U$ of $x_0$ in $X$ such that $B_n(U)\subset V$ for all $n\in\N$. But $\lim_n x_n=x_0$. Hence there is a $n_0\in\N$ such that $x_n\in U$ for all $n\ge n_0$. This implies that
\begin{equation*}
    B_n(x)= F(g_n,x_n)\in V\quad\hbox{for all}\quad n\ge n_0.
\end{equation*}
Thus $F$ is continuous in $(g_0,x_0)$.
\end{proof}

\begin{lem}\label{lem:H1_continuity}
The mapping
\begin{equation}\label{eq:Fcont}
    F\,:\, \D{2} \times \HH{1}\to\HH{1},\quad F(\varphi,v):=v\circ \varphi
\end{equation}
is continuous. Moreover, given $k,\,s\in\mathbb{N}$ with $s\ge 1$ and $k-s\ge 1$, the restriction of $F$ satisfies
\begin{equation*}
    F\in C(\D{k} \times \HH{k},\HH{k})\cap C(\D{k} \times \HH{k-s},\HH{k-s}).
\end{equation*}
\end{lem}

\begin{proof}
(a) By Sobolev's embedding theorem we know
that $\D{2}\hookrightarrow C^1(\Circle)$. Hence the chain rule ensures that $F$ is well-defined, i.e. $F(\varphi,v)\in\HH{1}$ for all $(\varphi,v)\in\D{2} \times \HH{1}.$ Moreover, fixing $\varphi\in\D{2}$, we have
\begin{equation*}
    F(\varphi,\cdot)\in\mathcal{L}(\HH{1},\HH{1}).
\end{equation*}

(b) Let now $v\in \HH{1}$ be fixed. We are going to show that
\begin{equation*}
    F(\cdot,v)\in C(\D{2},\HH{1}).
\end{equation*}
For this pick $\varphi_0\in\D{2}$ and $\varepsilon >0$. By Sobolev's embedding theorem, the function $v$ is uniformly continuous. Thus there is a $\delta>0$ such that
\begin{equation*}
    \abs{v(x)-v(y)} < \varepsilon \quad \text{for all}\quad \abs{x-y} <\delta.
\end{equation*}
Next let $j$ denote the embedding constant of $\HH{l}\hookrightarrow C(\Circle)$ for $l=1,\ 2$ and
choose $\varphi\in\D{2}$ such that $\norm{\varphi_0 - \varphi}_{H^2}<\delta/j$. Then
\begin{equation*}
    \abs{\varphi_0(x) - \varphi(x)} \le j \norm{\varphi_0 - \varphi}_{H^2} < \delta \quad \text{for all}\quad x\in\Circle.
\end{equation*}
Thus we get
\begin{equation}\label{eq:el2}
    \norm{v\circ\varphi_0 - v\circ\varphi}^2_{L^2} = \int_{\Circle} \abs{v(\varphi_0(x)) - v(\varphi(x))}^2 \, dx \le \varepsilon^2.
\end{equation}

To estimate $D(v\circ\varphi_0 -v\circ\varphi)$ in $L^2$, we first remark that it is
no restriction to assume that $\delta \in (0,1]$. Writing now $K:=j\Vert \varphi_0\Vert_{H^2}+1$ and
\begin{equation*}
    B_2(\delta) := \D{2}\cap\mathbb{B}_{H^2}(\varphi_0,\delta),
\end{equation*}
we have that
\begin{equation}\label{eq:xi-K}
    \norm{\varphi'}_{L^\infty} \le j\norm{\varphi}_{H^2} \le j\norm{\varphi_0}_{H^2} +1 = K \quad \text{for all} \quad \varphi\in B_2(\delta).
\end{equation}
Furthermore, letting $m(\varphi) := \norm{1/\varphi_x}_{L^\infty}$ for $\varphi\in\D{2}$, we have
\begin{equation*}
    \norm{f\circ\varphi}_{L^2}^2 \le m(\varphi)\norm{f}_{L^2}^2 \quad \text{for all} \quad f\in L^2(\Circle).
\end{equation*}
Note also that by shrinking $\delta >0$, we may assume that
\begin{equation}\label{eq:xi-sig}
    m(\varphi)\le 2 m(\varphi_0)\quad\hbox{for all}\quad \varphi\in B_2(\delta).
\end{equation}
We now proceed as follows. First we have
\begin{multline}\label{eq:triangle}
    \norm{\partial(v\circ\varphi_0 - v\circ\varphi)}^2_{L^2} \\
    \le \norm{v'\circ\varphi_0\cdot\varphi_0' - v'\circ\varphi_0\cdot\varphi'}^2_{L^2}
    + \norm{v'\circ\varphi_0\cdot\varphi' - v'\circ\varphi\cdot\varphi'}^2_{L^2}.
\end{multline}
For the first term of the right-hand side of \eqref{eq:triangle}, we find
\begin{equation}\label{eq:f-prim}
\begin{split}
    & \norm{v'\circ\varphi_0\cdot\varphi_0' - v'\circ\varphi_0\cdot\varphi'}^2_{L^2}
     = \int_{\Circle}\abs{v'(\varphi_0(x)}^2 \abs{\varphi_0'(x) - \varphi'(x)}^2 \, dx \\
    & \qquad \le\norm{\varphi_0' - \varphi'}^2_{L^\infty} \int_{\Circle}\abs{v'(\varphi_0(x))}^2 \, dx \\
    & \qquad \le j^2\,\Vert\varphi_0-\varphi\Vert^2_{H^2}\,m(\varphi_0)\,\Vert v\Vert^2_{H^1}=\delta^2\,m(\varphi_0)\,\Vert v\Vert^2_{H^1}.
\end{split}
\end{equation}

To estimate the second term in \eqref{eq:triangle}, choose $w\in C^2(\Circle)$ such that
\begin{equation} \label{eq:h1-K}
    \norm{v-w}_{H^1} \le \sqrt{\frac{1}{3\;m(\varphi_0)}} \; \frac{\varepsilon}{K} .
\end{equation}
Then we have
\begin{multline*}
    \norm{v'\circ\varphi_0 - v'\circ\varphi}_{L^2}^2 \le  \norm{v'\circ\varphi_0 - w'\circ\varphi_0}_{L^2}^2 + \\
    + \norm{w'\circ\varphi_0 - w'\circ\varphi}_{L^2}^2 + \norm{w'\circ\varphi - v'\circ\varphi}_{L^2}^2 \\
\begin{split}
   & \le (m(\varphi_0) + m(\varphi)) \norm{v'-w'}_{L^2}^2 + \int_{\Circle}\abs{w'(\varphi_0(x)) - w'(\varphi(x))}^2 dx \\
   & \le 3m(\varphi_0) \norm{v'-w'}_{L^2}^2 + \norm{w''}_{L^\infty}^2 \int_{\Circle} \abs{\varphi_0(x) - \varphi(x)}^2 dx,
\end{split}
\end{multline*}
where we also employed the mean value theorem and \eqref{eq:xi-sig} to derive the last estimate. Invoking \eqref{eq:xi-K} and \eqref{eq:h1-K}, we get
\begin{equation}\label{eq:penultimo}
    \norm{v'\circ\varphi_0\cdot\varphi' - v'\circ\varphi\cdot\varphi'}_{L^2} \le \varepsilon + \delta K \norm{w''}_{L^\infty}
\end{equation}
for all $\varphi\in B_2(\delta)$. Combining \eqref{eq:el2}, \eqref{eq:triangle}, \eqref{eq:f-prim}, and \eqref{eq:penultimo}, we arrive at the following estimate
\begin{equation}\label{eq:ultimo}
    \norm{v\circ\varphi_0 - v\circ\varphi}_{H^1} \le 2\,\varepsilon + \delta \left(\sqrt{m(\varphi_0)} \norm{v}_{H^1}
    + K\,\norm{w''}_{L^\infty}\right)
\end{equation}
for all $\varphi\in B_2(\delta)$. Shrinking $\delta>0$, we get from \eqref{eq:ultimo} that
\begin{equation*}
    \norm{v\circ\varphi_0 - v\circ\varphi}_{H^1} \le 3\,\varepsilon
\end{equation*}
for all $\varphi\in B_2(\delta)$. Thus $F(\cdot,v)$ is continuous in $\varphi_0\in\D{2}$. Invoking Lemma \ref{lem:sep_continuity}, we find that $F\in C(\D{2}\times\HH{1},\HH{1})$.

(c) Let $k\ge 2$ be given. Then it follows from the considerations from \cite[page 108]{EM70} that
\begin{equation*}
    F(\varphi,\cdot)\in\mathcal{L}(\HH{k},\HH{k}),
\end{equation*}
for all $\varphi\in\D{k}$ and that
\begin{equation*}
    F(\cdot,v)\in C(\D{k},\HH{k}),
\end{equation*}
for all $v\in\HH{k}$. Hence, again by lemma~\ref{lem:sep_continuity}, we conclude that
\begin{equation*}
    F\in C(\D{k}\times\HH{k},\HH{k}).
\end{equation*}
The last assertion is now obvious.
\end{proof}

\begin{rem}
(a) For simplicity we treated here the case $s\in\mathbb{N}$. Using an intrinsic representation of the Sobolev norm for $s\in\mathbb{R}$ with  $s\ge 1$, it is possible to extend the results of Lemma \ref{lem:H1_continuity} to non-integer values of $s\ge 1$.

(b) A similar result to \eqref{eq:Fcont} has recently been established in \cite{DKT07}. However, on the one hand, Corollary 3 in \cite{DKT07} fits not precisely into our setting, and on the other hand our scale of Sobolev spaces is simpler than the one in \cite{DKT07}. Therefore we decided to present a self-contained proof of \eqref{eq:Fcont}.

(c) The higher the spatial regularity in the group $\D{k}$ and the Lie algebra $\HH{k}$, the better the regularity of the mapping $F$ in
lemma~\ref{lem:H1_continuity}, cf. \cite{EM70}. However, we are not aware of better regularity of $F$ than \eqref{eq:Fcont}. Finally, we remark that  the continuity of $F$ is sufficient for our purposes.
\end{rem}

% ----------------------------------------------------------------

\section{Proof of Proposition~\ref{prop:smoothness_of_lambda}}
\label{sec:proof_prop}

In this section we provide the completion of the proof of the smoothness of the inertia operator $\Lambda_\varphi(v)$ with respect to suitable Sobolev norms.

\begin{lem}\label{lem:nth_derivative_symbol}
Let $P$ be a Fourier multiplier on $\CS$, and let $P_{n}$ be the multilinear operator defined in Proposition~\eqref{lem:nth_derivative}
for some $n\in\mathbb{N}$. Then we have
\begin{equation}\label{eq:nth_derivative_symbol}
    P_{n}(\e_{m_{1}}, \dotsc ,\e_{m_{n}})\e_{m_{0}} = p_{n}(m_{0}, m_{1}, \dotsc , m_{n}) \e_{m_{0} + m_{1} \dotsb + m_{n}},
\end{equation}
where the sequence $p_{n}$ is defined inductively by $p_{0} = p$ (the symbol of $P$) and
\begin{multline}\label{eq:derivative_symbol_recurrence}
    p_{n+1}(m_{0}, \dotsc , m_{n+1}) = (2 \pi i) \Big[ (m_{0} + \dotsb + m_{n}) p_{n}(m_{0}, \dotsc , m_{n}) \\
    - \sum_{j=0}^{n} m_{j} \, p_{n}(m_{0}, \dotsc ,  m_{j} + m_{n+1}, \dotsc ,  m_{n})\Big],
\end{multline}
and $m_j\in\mathbb{Z}\setminus\{0\}$, $j=1,\dots,n.$
\end{lem}

\begin{rem}
For $P := \Lambda = H \circ D$, we have
\begin{equation}\label{eq:principal_symbol}
    p_{0}(m_{0}) = \abs{m_{0}},
\end{equation}
and
\begin{equation}\label{eq:first_derivative_symbol}
    p_{1}(m_{0}, m_{1}) = (2 \pi i) m_{0} \Big( \abs{m_{0}} - \abs{m_{0} + m_{1}} \Big)
\end{equation}
and
\begin{multline}\label{eq:second_derivative_symbol}
    p_{2}(m_{0}, m_{1}, m_{2}) = (2 \pi i)^{2} m_{0} \Big( (m_{0} + m_{1} + m_{2}) \abs{m_{0} + m_{1} + m_{2}} \\
        - (m_{0} + m_{1}) \abs{m_{0} + m_{1}} - (m_{0} + m_{2}) \abs{m_{0} + m_{2}} + m_{0} \abs{m_{0}} \Big).
\end{multline}
\end{rem}

\begin{proof}
Invoking Lemma \ref{lem:FOp}, the case  $n=0$ is clear.
Suppose that equation~\eqref{eq:nth_derivative_symbol}
is true
for some $n \ge 0$. Then, using recurrence relation~\eqref{eq:recurrence_relation}, we have
\begin{multline*}
    P_{n+1}(\e_{m_{1}}, \dotsc , \e_{m_{n+1}})\e_{m_{0}} = \e_{m_{n+1}}D\big( P_{n}(\e_{m_{1}}, \dotsc , \e_{m_{n}})\e_{m_{0}} \big) \\
    - P_{n}(\e_{m_{1}}, \dotsc , \e_{m_{n}}) \big( \e_{m_{n+1}}D \e_{m_{0}}\big) - \sum_{j=1}^{n}P_{n}(\e_{m_{1}}, \dotsc ,D\e_{m_{j}}\e_{m_{n+1}}, \dotsc , \e_{m_{n}}),
\end{multline*}
which is equal to
\begin{multline*}
    (2 \pi i) \Big\{ (m_{0} + \dotsb + m_{n}) p_{n}(m_{0}, \dotsc , m_{n}) - m_{0} \, p_{n}(m_{0} + m_{n+1}, \dotsc , m_{n}) \\
    - \sum_{j=1}^{n} m_{j}\, p_{n}(m_{0}, \dotsc , m_{j} + m_{n+1}, \dotsc ,m_{n})\Big\}\e_{m_{0} + \dotsb + m_{n+1}}.
\end{multline*}
This shows that equation~\eqref{eq:nth_derivative_symbol} is true for $n+1$ with
\begin{multline*}
    p_{n+1}(m_{0}, \dotsc , m_{n+1}) = (2 \pi i) \Big[ (m_{0} + \dotsb + m_{n}) p_{n}(m_{0}, \dotsc , m_{n}) \\
    - \sum_{j=0}^{n} m_{j} \, p_{n}(m_{0}, \dotsc ,  m_{j} + m_{n+1}, \dotsc ,  m_{n})\Big]
\end{multline*}
and achieves the proof.
\end{proof}

\begin{lem}\label{lem:nth_derivative_estimates}
Let $P$ be a Fourier multiplier of order $s \in \N$ and $k \ge s + 1$. Let $P_{n}$ be the $(n+1)$-multilinear operator defined by the recurrence relation~\eqref{eq:recurrence_relation} with $P_{0}=P$. Suppose that there exists a constant $C_{n} >0 $, such that
\begin{equation}\label{eq:nth_derivative_estimates}
    \abs{p_{n}(m_{0}, \dotsc , m_{n})} \le C_{n} \abs{m_{0}}^{s} \dotsb \abs{m_{n}}^{s}
\end{equation}
for all $m_{j} \in \Z \setminus\{0\}$. Then $P_{n}$ extends to a bounded multilinear operator
\begin{equation*}
    P_{n} : \overbrace{\HH{k} \times \dotsb \times \HH{k}}^{n+1} \to H^{k-s}(\Circle).
\end{equation*}
\end{lem}

\begin{proof}
By virtue of Proposition~\ref{lem:nth_derivative_symbol}, we have
\begin{multline*}
    \norm{P_{n}(u_{1}, \dotsc ,u_{n})u_{0}}_{H^{k-s}}^{2} = \\
    \sum_{l\in \Z} \abs{\sum_{m_{0} + \dotsb + m_{n} = l} \hat{u}_{0}(m_{0}) \dotsb \hat{u}_{n}(m_{n}) p_{n}(m_{0}, \dotsc , m_{n})}^{2}\norm{\e_{l}}_{H^{k-s}}^{2},
\end{multline*}
for any smooth functions $u_{0}, u_{1}, \dotsc , u_{n}$, since $(\e_{l})_{l \in \Z}$ is an orthogonal system for the $H^{k-s}$ inner product. Therefore, if relation~\eqref{eq:nth_derivative_estimates} is satisfied, we get
\begin{multline*}
    \norm{P_{n}(u_{1}, \dotsc ,u_{n})u_{0}}_{H^{k-s}}^{2} \le \\
    C_{n}\sum_{l\in \Z} \left(\sum_{m_{0} + \dotsb + m_{n} = l} \abs{m_{0}}^{s}\abs{\hat{u}_{0}(m_{0})} \dotsb \abs{m_{n}}^{s}\abs{\hat{u}_{n}(m_{n})} \right)^{2}\norm{\e_{l}}_{H^{k-s}}^{2}.
\end{multline*}
Observe now that, given smooth functions $v_{0}, v_{1}, \dotsc , v_{n}$, we have
\begin{equation*}
    \widehat{v_{0} \dotsb v_{n}}(l) = \sum_{m_{0} + \dotsb + m_{n} = l} \hat{v}_{0}(m_{0}) \dotsb \hat{v}_{n}(m_{n}).
\end{equation*}
In addition $\HH{k-s}$ is a Banach algebra, since $k-s \ge 1$. Consequently there exists a constant $C^{\prime}_{n,k,s}$ such that
\begin{multline*}
    \sum_{l\in \Z} \abs{\sum_{m_{0} + \dotsb + m_{n} = l} \hat{v}_{0}(m_{0}) \dotsb \hat{v}_{n}(m_{n}) }^{2}\norm{e_{l}}_{H^{k-s}}^{2} \\
    \le C^{\prime}_{n,k,s}\norm{v_{0}}_{H^{k-s}}^{2} \dotsb \norm{v_{n}}_{H^{k-s}}^{2}
\end{multline*}
for every smooth functions $v_{0}, v_{1}, \dotsc , v_{n}$.
Putting now $\hat{v}_{p}(m_{p}) = \abs{m_{p}^{s}\hat{u}_{p}(m_{p})}$ in this last inequality and using the fact that the functions with Fourier coefficient $\hat{u}(m)$ and $\abs{\hat{u}(m)}$ have the same $H^{k-s}$ norm, we obtain
\begin{align*}
    \norm{P_{n}(u_{1}, \dotsc ,u_{n})u_{0}}_{H^{k-s}}^{2} & \le C_{n}C^{\prime}_{n,k,s}\norm{u_{0}^{(s)}}_{H^{k-s}}^{2} \dotsb \norm{u_{n}^{(s)}}_{H^{k-s}}^{2} \\
    & \le C^{\prime\prime}_{n,k,s}\norm{u_{0}}_{H^{k}}^{2} \dotsb \norm{u_{n}}_{H^{k}}^{2},
\end{align*}
which achieves the proof.
\end{proof}

\begin{cor}\label{cor:continuity_derivatives}
Let $P$ be a Fourier multiplier of order $s$. Let $r\in\N$ and $k \ge s + 1$.
Suppose that the operators $P_{n}$, defined in Proposition~\eqref{lem:nth_derivative}, extend to \emph{bounded} multilinear operators
\begin{equation*}
    P_{n} : \overbrace{\HH{k} \times \dotsb \times \HH{k}}^{n+1} \to H^{k-s}(\Circle).
\end{equation*}
for $0 \le n \le r$. Then
\begin{equation*}
	(\varphi, v) \mapsto P_{\varphi} (v) = R_{\varphi} \circ P \circ R_{\varphi^{-1}} (v).
\end{equation*}
is of class $C^{r}$ from $\D{k} \times \HH{k}$ to $\HH{k-s}$.
\end{cor}

\begin{proof}
Notice first that if $P_{0} = P$ is bounded, then $(\varphi, v) \mapsto P_{\varphi} (v)$ is continuous from $\D{k} \times \HH{k}$ into $\HH{k-s}$, by virtue of lemma~\ref{lem:H1_continuity}. Suppose now that the bilinear operator $P_{1}$ is bounded and let
\begin{equation*}
   P_{\varphi}^{1}(u,v) := \big(R_{\varphi} P_{1}(u)R_{\varphi^{-1}} \big)(v).
\end{equation*}
Applying Lemma~\ref{lem:H1_continuity}, we deduce that the map
\begin{equation*}
    \D{k} \times \HH{k} \times \HH{k} \to \HH{k-s}, \qquad (\varphi, u, v) \mapsto P_{\varphi}^{1}(u,v)
\end{equation*}
is continuous. Using Proposition \ref{lem:nth_derivative} and the mean value theorem in the global chart $U$ defined by \eqref{eq:smooth_chart}, we get
\begin{equation}\label{eq:mean_value}
   P_{\varphi +u}(v)  - P_{\varphi}(v)  = \int^{1}_{0} P^{1}_{\varphi +tu}(u,v) \, dt
\end{equation}
for smooth maps $\varphi, u, v$. But, since both sides of \eqref{eq:mean_value} are continuous in all the variables, we deduce, using a density argument, that this relation is still true for $\varphi \in \D{k}$, $u,v\in \HH{k}$. We conclude therefore that $(\varphi,v) \mapsto  P_{\varphi}(v)$ is a $C^{1}$ map from $\D{k} \times \HH{k}$ into $\HH{k-s}$. An inductive argument using the same reasoning for $P_{n}$ shows that $(\varphi,v) \mapsto  P_{\varphi}(v)$ is a $C^{n}$ map from $\D{k} \times \HH{k}$ into $\HH{k-s}$, for each $n \le r$.
\end{proof}

Finally, we need the following elementary lemma.

\begin{lem}\label{lem:n_order_estimate}
Let $f : \R \to \R$ be a function of class $C^{n-1}$ and such that $f^{(n-1)}$ satisfies a Lipschitz condition with Lipschitz constant $K$. Then
\begin{equation*}
    \abs{\sum_{p=0}^{n} (-1)^{p} \sum_{\substack{I \subset \set{1, \dotsc , n}, \\ \abs{I} = p}} f\big(t + \sum_{j \in I} m_{j}\big)} \le K \prod_{j=1}^{n}\abs{m_{j}},
\end{equation*}
for all $t \in \R$ and all $m_{j}\in \R$.
\end{lem}

\begin{proof}
Let $g_{k}$ be the sequence of functions defined inductively by
\begin{equation*}
    g_{1}(t) = f(t+m_{1}) - f(t), \qquad g_{k+1}(t) = g_{k}(t+m_{k}) - g_{k}(t).
\end{equation*}
Then, we have
\begin{equation*}
    g_{n}(t) = (-1)^{n}\sum_{p=0}^{n} (-1)^{p} \sum_{\substack{I \subset \set{1, \dotsc , n}, \\ \abs{I} = p}} f\big(t + \sum_{j \in I} m_{j}\big)
\end{equation*}
On the other hand, the Lipschitz condition on the $(n-1)$ derivative of $f$ leads to
\begin{equation*}
    \abs{g_{1}^{(n-1)}(t)} \le K \abs{m_{1}}, \qquad \forall t \in \R .
\end{equation*}
Now, using inductively the mean value theorem, we get
\begin{equation*}
    \abs{g_{k}^{(n-k)}(t)} \le K \abs{m_{1}} \dotsb \abs{m_{k}}, \qquad \forall t \in \R .
\end{equation*}
In particular, for $k = n$, we have
\begin{equation*}
    \abs{g_{n}(t)} \le K \prod_{j=1}^{n}\abs{m_{j}}, \qquad \forall t \in \R ,
\end{equation*}
which achieves the proof.
\end{proof}

\begin{proof}[Proof of Proposition~\ref{prop:smoothness_of_lambda}]
For each $n \ge 1$, let $f_{n}(t) = t^{n-1}\abs{t}$. Then $f_{n}$ is of class $C^{n-1}$ on $\R$ and $f_{n}^{(n-1)}$ satisfies a global Lipschitz condition with Lipschitz constant $(n-1)!$. We are going to show that
\begin{equation}\label{eq:main_equality}
    p_{n}(m_{0}, m_{1}, \dotsc , m_{n}) = (2 \pi i )^{n} m_{0} \sum_{p=0}^{n} (-1)^{p} \sum_{\substack{I \subset \set{1, \dotsc , n}, \\ \abs{I} = p}} f_{n}\big(m_{0} + \sum_{j \in I} m_{j}\big),
\end{equation}
for each $n \ge 1$. Then, by virtue of Corollary~\ref{cor:continuity_derivatives} and Lemma~\ref{lem:n_order_estimate}, this will demonstrate that
\begin{equation*}
	(\varphi, v) \mapsto \Lambda_{\varphi} (v) = R_{\varphi} \circ \Lambda \circ R_{\varphi^{-1}} (v)
\end{equation*}
is smooth from $\D{k} \times \HH{k}$ to $\HH{k-1}$.

For $n=1$, we have
\begin{equation*}
    p_{1}(m_{0}, m_{1}) = (2 \pi i ) m_{0} \big( \abs{m_{0}} - \abs{m_{0} + m_{1}} \big)
\end{equation*}
so equation~(\ref{eq:main_equality}) is true for $n=1$. Now, suppose inductively that this equation is valid for some $n\ge 1$. Using the recurrence relation~(\ref{eq:derivative_symbol_recurrence}), we get
\begin{multline*}
    p_{n+1}(m_{0}, m_{1}, \dotsc , m_{n+1}) = (2 \pi i )^{n+1} m_{0} \sum_{p=0}^{n} (-1)^{p} \sum_{\substack{I \subset \set{1, \dotsc , n}, \\ \abs{I} = p}}  \Big\{ \\
    (m_{0} + \dotsb + m_{n}) f_{n}\big(m_{0} + \sum_{j \in I} m_{j}\big) - \sum_{k=1}^{n} m_{k} \, f_{n}\big(m_{0} + \sum_{j \in I} m_{j} + \delta_{I}(k) \, m_{n+1} \big)  \\
    - (m_{0} + m_{n+1})  f_{n}\big(m_{0} + \sum_{j \in I} m_{j} + m_{n+1} \big) \Big\} ,
\end{multline*}
which can be rewritten as
\begin{multline*}
    (2 \pi i )^{n+1} m_{0} \sum_{p=0}^{n} (-1)^{p} \sum_{\substack{I \subset \set{1, \dotsc , n}, \\ \abs{I} = p}}  \Big\{ \big( m_{0} + \sum_{j \in I} m_{j} \big) f_{n}\big(m_{0} + \sum_{j \in I} m_{j}\big) \\
    - \big( m_{0} + \sum_{j \in I} m_{j} + m_{n+1} \big) \, f_{n}\big(m_{0} + \sum_{j \in I} m_{j} + m_{n+1} \big) \Big\} .
\end{multline*}
using the fact that $f_{n+1}(t) = tf_{n}(t)$, we have therefore
\begin{multline*}
    p_{n+1}(m_{0}, m_{1}, \dotsc , m_{n+1}) = (2 \pi i )^{n+1} m_{0} \sum_{p=0}^{n} (-1)^{p} \sum_{\substack{I \subset \set{1, \dotsc , n}, \\ \abs{I} = p}}  \Big\{ \\
    f_{n+1}\big(m_{0} + \sum_{j \in I} m_{j}\big) - f_{n+1}\big(m_{0} + \sum_{j \in I} m_{j} + m_{n+1} \big) \Big\} ,
\end{multline*}
which is equal to
\begin{multline*}
    (2 \pi i )^{n+1} m_{0} \Big\{ \sum_{p=0}^{n} (-1)^{p} \sum_{\substack{I \subset \set{1, \dotsc , n+1}, \\ \abs{I} = p,\, n+1 \notin I} }
    f_{n+1}\big(m_{0} + \sum_{j \in I} m_{j}\big) \\
    +  \sum_{p=0}^{n} (-1)^{p+1} \sum_{\substack{I \subset \set{1, \dotsc , n+1}, \\ \abs{I} = p+1,\, n+1 \in I} } f_{n+1} \big(m_{0} + \sum_{j \in I} m_{j}\big) \Big\} .
\end{multline*}
But this last expression is exactly
\begin{equation*}
    (2 \pi i )^{n+1} m_{0} \sum_{p=0}^{n+1} (-1)^{p} \sum_{\substack{I \subset \set{1, \dotsc , n+1}, \\ \abs{I} = p}} f_{n+1}\big(m_{0} + \sum_{j \in I} m_{j}\big),
\end{equation*}
which achieves the proof.
\end{proof}

% ----------------------------------------------------------------

\section*{Acknowledgments}

The authors wish to express their gratitude to the Erwin Schr\"{o}dinger International Institute for Mathematical Physics for providing an excellent research environment during the program ``Integrable systems of hydrodynamic type'' (Oct. 12--23, 2009). M. W. acknowledges financial support by the JSPS Postdoctoral Fellowship P09024.

Finally it is a pleasure to thank Anders Melin and Elmar Schrohe for helpful discussions.

% ----------------------------------------------------------------

\end{document}